\documentclass[12pt,a4paper]{amsart}

\usepackage{amssymb}
\usepackage{calrsfs}
\usepackage{url}
\usepackage{multirow}
\newcommand{\QQ}{{\mathbb{Q}}}
\newcommand{\FF}{{\mathbb{F}}}

\newcommand{\RR}{{\mathbb{R}}}
\newcommand{\ZZ}{{\mathbb{Z}}}

\newcommand{\fB}{{\mathfrak{B}}}

\newcommand{\cG}{{\mathcal{G}}}
\newcommand{\cH}{{\mathcal{H}}}
\newcommand{\cM}{{\mathcal{M}}}
\newcommand{\cW}{{\mathcal{W}}}

\newcommand{\cN}{{\mathcal{N}}}
\newcommand{\cE}{{\mathcal{E}}}
\newcommand{\cI}{{\mathcal{I}}}
\newcommand{\cS}{{\mathcal{S}}}

\newcommand\Irr{\operatorname{Irr}}
\newcommand\GL{\operatorname{GL}}
\newcommand\Uch{\operatorname{Uch}}

\newtheorem{thm}{Theorem}[section]

\newtheorem{lem}[thm]{Lemma}

\theoremstyle{definition}

\newtheorem{exmp}[thm]{Example}

\theoremstyle{remark}
\newtheorem{rem}[thm]{Remark}

\numberwithin{equation}{section}
\renewcommand{\leq}{\leqslant}
\renewcommand{\geq}{\geqslant}

\begin{document}

%%%%%%%%%%%%%%%%%%%%%%%%%%%%%%%%%%%%%%%%%%%%%%%%%%%%%%%%%%%%%%%%%%
\title{On the labelling of characters of Weyl groups of type $F_4$}
%%%%%%%%%%%%%%%%%%%%%%%%%%%%%%%%%%%%%%%%%%%%%%%%%%%%%%%%%%%%%%%%%%

\author{Meinolf Geck and Jonas Hetz}
\address{FB Mathematik, Universit\"at Stuttgart, 
Pfaffenwaldring 57, 70569 Stuttgart, Germany.}
\email{meinolf.geck@mathematik.uni-stuttgart.de}
\address{Lehrstuhl für Algebra und Zahlentheorie, RWTH Aachen, 
Pontdriesch 14/16, 52062 Aachen, Germany.}
\email{jonas.hetz@rwth-aachen.de}

%\thanks{The authors gratefully acknowledge support by the Deutsche
% Forschungsgemeinschaft (DFG, German Research Foundation) --- Project-ID
% 286237555 -- TRR 195.}

\begin{abstract} In the literature on finite groups of Lie type, there 
exist two different conventions about the labelling of the irreducible
characters of Weyl groups of type~$F_4$. We point out some issues 
concerning these two conventions and their effect on tables about 
unipotent characters or the Springer correspondence. Using experiments 
related to these issues with the computer algebra system {\sf CHEVIE},
we spotted an error in Spaltenstein's tables for the generalised Springer
correspondence in type~$E_7$.
\end{abstract}

%\keywords{Weyl groups, character theory, groups of Lie type}


%\date{\today}

\maketitle

\pagestyle{myheadings}
\markboth{On the labelling of characters of type $F_4$}{On the 
labelling of characters of type $F_4$}

%%%%%%%%%%%%%%%%%%%%%%%%%%%%%%%%%%%%%%%%%%%%%%%%%%%%%%%%%%%%%%%%%%%%%%%%%%%
\section{Introduction} \label{sec0}

This note is about certain notational conventions in the representation 
theory of finite groups of Lie type. General references are the books by 
Carter \cite{Ca2}, Digne--Michel \cite{rDiMi} and Lusztig \cite{LuB};
see also \cite{gema} for a more recent survey. There are some issues 
in relation to Weyl groups of type $F_4$. These groups (and related 
data like generic degrees, Springer representations etc.) occur in a 
substantial way in a number of situations including:
\begin{itemize}
\item The parametrisation of unipotent characters of a Chevalley group 
$F_4(q)$.
\item Degree formulae for unipotent principal series characters of a twisted
Chevalley group ${^2\!E}_6(q)$.
\item Degrees of unipotent characters of $E_8(q)$ in the Harish-Chandra 
series above a cuspidal unipotent character of a Levi subgroup of type
$D_4$. 
\item The Springer correspondence for a simple algebraic group of type 
$F_4$.
\item The generalised Springer correspondence for a simple algebraic group 
of type $E_7$ with respect to a cuspidal datum in a Levi subgroup 
of type $A_1{\times} A_1{\times} A_1$.
\item The generalised Springer correspondence for a simple algebraic group 
of type $E_8$ with respect to a cuspidal datum in a Levi subgroup 
of type $D_4$.
\end{itemize}
In a root sytem of type $F_4$ there are long roots and short roots, but the
picture is highly symmetric and, hence, it often requires some extra 
care to fix notation and conventions. Unfortunately, this is not always 
done consistently in the existing literature, which prompted us to write 
this note. In fact, during the work on \cite{Het3}, there actually occurred
some delicate contradictions as a result of a misunderstanding regarding
the conventions used in existing tables. (This will be briefly explained
in Example~\ref{gspr2}.). Our purpose here is to try to clarify some points 
in relation to the long/short root issue in type $F_4$. In the course
of verifying a number of tables, we actually discovered an error in
Spaltenstein's table for the generalised Springer correspondence in type
$E_7$; see Lemma~\ref{esp}.

%%%%%%%%%%%%%%%%%%%%%%%%%%%%%%%%%%%%%%%%%%%%%%%%%%%%%%%%%%%%%%%%%%%%%%%%%%%
\section{Characters of a Coxeter group of type $F_4$} \label{sec1}

Let $W$ be a Coxeter group of type $F_4$. The basic and widely used reference
for the character table of $W$ is Kondo \cite{kondo}. In that article, $W$
is constructed in a purely group-theoretical way, without reference to an
underlying root system. With the notation of \cite[Prop.~5.1]{kondo},
we have $W=\langle d,a, \tau,\tau\sigma\rangle$ where $d,a,\tau,\tau\sigma$
are elements of order $2$ that satisfy the following braid relations:
\[ (d{\cdot} a)^3=(\tau{\cdot}\tau\sigma)^3=(a{\cdot}\tau)^4=1, \quad 
d{\cdot} \tau=\tau {\cdot} d,\quad d{\cdot}\tau\sigma=\tau\sigma{\cdot} d,
\quad a{\cdot} \tau \sigma=\tau\sigma{\cdot} a.\]
Thus, the corresponding Coxeter diagram is as follows.
\begin{center}
\begin{picture}(220,20)
\put( 62,12){$d$}
\put( 92,12){$a$}
\put(122,12){$\tau$}
\put(149,12){$\tau\sigma$}
\put( 65,04){\circle*{6}}
\put( 95,04){\circle*{6}}
\put(125,04){\circle*{6}}
\put(155,04){\circle*{6}}
\put( 65,04){\line(1,0){30}}
\put( 95,04){\line(1,0){30}}
\put(125,04){\line(1,0){30}}
\put(107,08){\scriptsize{$4$}}
\end{picture}
\end{center}
The group $W$ has $25$ irreducible characters, of degrees $1,2,4,6,8,9,12,
16$; we shall denote by $n_j$ the $j$-th character of degree $n$ in
Kondo's table \cite[p.~152]{kondo}. (This is a standard notation in the
existing literature, e.g., Alvis \cite{alvis}, Carter \cite{Ca2}, Lusztig
\cite{LuB}; note that Shoji \cite{Sho}, Spaltenstein \cite{Spa1} use
a slightly different notation for the characters of degree~$4$.)

\begin{table}[htbp]\caption{Irreducible characters for type $F_4$}
\label{tabf4}
\begin{center} $\begin{array}{ccrrrr} \hline \phi&&
\;\;a_\phi&\;\;d,a&\;\;\tau,\tau\sigma & \;\;a\tau\\ \hline 
1_1 & \phi_{1,0}  & 0 &  1 & 1 & 1\\ \hline 
4_2 & \phi_{4,1}  & 1 &  2 & 2 & 2\\ 
2_1 & \phi_{2,4}''  & 1 &  2 & 0 & 0\\ 
2_3 & \phi_{2,4}'  & 1 &  0 & 2 & 0 \\ \hline 
9_1 & \phi_{9,2}  & 2 &  3 & 3 & 1\\ \hline 
8_1 & \phi_{8,3}''  & 3 &  4 & 0 & 0\\ \hline
8_3 & \phi_{8,3}'  & 3 &  0 & 4 & 0\\ \hline 
12_1 & \phi_{12,4}  & 4 &  0 & 0 & 0\\ 
16_1 & \phi_{16,5}  & 4 &  0 & 0 & 0 \\ 
9_2 & \phi_{9,6}''  & 4 &  3 & -3 & -1\\ 
6_2 & \phi_{6,6}''  & 4 &  0 & 0 & 2\\ 
9_3 & \phi_{9,6}'  & 4 &  -3 & 3 & -1\\ 
6_1 & \phi_{6,6}'  & 4 &  0 & 0 & -2 \\ 
\hline \end{array}\qquad\qquad\begin{array}{ccrrrr} \hline \phi&&
\;\;a_\phi&\;\;d,a&\;\;\tau,\tau\sigma & a\tau\\ \hline 
4_4 & \phi_{4,7}'  & 4 &  -2 & 2 & -2 \\ 
4_3 & \phi_{4,7}''  & 4 &  2 & -2  & -2\\ 
4_1 & \phi_{4,8}  & 4 &  0 & 0  & 0\\ 
1_3 & \phi_{1,12}'  & 4 &  -1 & 1  & -1\\ 
1_2 & \phi_{1,12}''  & 4 &  1 & -1 & -1\\ \hline 
8_2 & \phi_{8,9}'  & 9 &  -4 & 0  & 0 \\  \hline
8_4 & \phi_{8,9}''  & 9 &  0 & -4  & 0 \\ \hline 
9_4 & \phi_{9,10}  & 10 &  -3 & -3 & 1\\ \hline 
4_5 & \phi_{4,13}  & 13 &  -2 & -2 & 2\\ 
2_2 & \phi_{2,16}'  & 13 &  -2 & 0 & 0 \\ 
2_4 & \phi_{2,16}''  & 13 &  0 & -2 & 0 \\ \hline 
1_4 & \phi_{1,24}  & 24 &  -1 & -1  & 1\\ 
\hline & & & & \end{array}$\end{center}
\end{table}

Table~\ref{tabf4} contains some information about the values of the
characters on $d,a,\tau,\tau\sigma$ and $a\tau$. The second column of that 
table refers to the notation introduced by Carter \cite[p.~413]{Ca2} (where
one can also find the complete character table); the third column refers
to the $a$-function, as defined in \cite[4.1]{LuB}, and printed in
\cite[4.10]{LuB}.

\begin{rem} \label{rem12} There is a group automorphism $\iota\colon 
W\rightarrow W$ such that $\iota(d)=\tau\sigma$, $\iota(a)=\tau$, $\iota(\tau)
=a$ and $\iota(\tau\sigma)=d$. This automorphism has order $2$ and it induces
a permutation of the irreducible characters of $W$ as follows:
\[ 1_2\leftrightarrow 1_3,\quad 2_1\leftrightarrow 2_3,\quad
2_2\leftrightarrow 2_4,\quad 4_3\leftrightarrow 4_4,\quad
8_1\leftrightarrow 8_3,\quad 8_2\leftrightarrow 8_4,\quad
9_2\leftrightarrow 9_3;\]
or, with Carter's notation:
\begin{gather*}
\phi_{1,12}'\leftrightarrow \phi_{1,12}'',\quad
 \phi_{2,4}'\leftrightarrow \phi_{2,4}'',\quad
 \phi_{2,16}'\leftrightarrow \phi_{2,16}'',\quad
 \phi_{4,7}'\leftrightarrow \phi_{4,7}'',\\\quad
 \phi_{8,3}'\leftrightarrow \phi_{8,3}'',\quad
 \phi_{8,9}'\leftrightarrow \phi_{8,9}'',\quad
 \phi_{9,6}'\leftrightarrow \phi_{9,6}''.
\end{gather*}
(All other characters are fixed.)
\end{rem}

\begin{rem} \label{rem12a} Let $(\phi',\phi'')$ be a pair of characters
that are permuted by the above automorphism $\iota\colon W \rightarrow W$.
By inspection of Table~\ref{tabf4} we see that we always have $\phi'(a)=
\phi'(d)\leq 0$ and $\phi''(a)=\phi''(d)\geq 0$. (Hence, this property
characterises $\phi'$ and~$\phi''$.)
\end{rem}

\begin{rem} In the following sections, the Coxeter group $W$ will typically
arise as the Weyl group of a (crystallographic) root system $\Phi$ of type
$F_4$ in some finite-dimensional Euclidean vector space~$E$ with scalar 
product $(\;,\;)\colon E\times E\rightarrow \RR$. Such a root system 
contains roots of different lengths, which gives rise to the long/short root 
issue mentioned in the introduction. Let us fix a set of simple roots 
$\Pi=\{\alpha_1,\alpha_2,\alpha_3,\alpha_4\}\subseteq \Phi$, labelled 
in such a way that $\alpha_1,\alpha_2$ are long and $\alpha_3,\alpha_4$ are 
short:
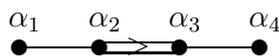
\begin{center}
\begin{picture}(220,20)
\put( 61,12){$\alpha_1$}
\put( 91,12){$\alpha_2$}
\put(121,12){$\alpha_3$}
\put(151,12){$\alpha_4$}
\put( 65,04){\circle*{6}}
\put( 95,04){\circle*{6}}
\put(125,04){\circle*{6}}
\put(155,04){\circle*{6}}
\put( 65,04){\line(1,0){30}}
\put( 95,02){\line(1,0){30}}
\put( 95,06){\line(1,0){30}}
\put(125,04){\line(1,0){30}}
\put(105,1.5){$>$}
\end{picture}
\end{center}
\smallskip
Then, denoting by $s_i\in\GL(E)$ the reflection with root $\alpha_i$ 
($1\leqslant i\leqslant4$), we have $W\cong\langle s_1,s_2,s_3,s_4\rangle
\subseteq\GL(E)$. Hence, when referring to the labelling of $\Irr(W)$ in 
Table~\ref{tabf4}, it is necessary to match the above generators $s_1,s_2,
s_3,s_4$ to Kondo's generators $d,a,\tau,\tau\sigma$; in other words, one 
has to specify which of Kondo's generators should correspond to long roots, 
and which to short roots. There are precisely two possibilities:

\smallskip
\begin{center}
\fbox{\begin{tabular}{r@{\hspace{1pt}}l} (AC): &  $\;\;s_1{=}d$, $s_2{=}a$, 
$s_3{=}\tau$, $s_4{=}\tau\sigma, \;$ that is, $d,a$ are reflections in long 
roots;\\ or (L): &  $\;\;s_1{=}\tau\sigma$, $s_2{=}\tau$, $s_3{=}a$, 
$s_4{=}d,\;$ that is, $d,a$ are reflections in short roots.\end{tabular}}
\end{center} 

\smallskip
Both (AC) and (L) exist in the literature: In \cite[4.10]{LuB}, Lusztig 
explicitly identifies $d,a$ with reflections in long simple coroots and 
$\tau,\tau\sigma$ with reflections in short simple coroots. Recall that 
the coroots are defined by $\alpha^\vee:=2\alpha/(\alpha,\alpha)\in E$ for 
$\alpha\in \Phi$. Then $\Phi^\vee=\{\alpha^\vee\mid \alpha\in \Phi\}$ also
is a root system of type $F_4$ in $E$, with simple roots given by $\Pi^\vee=
\{\alpha_1^\vee,\alpha_2^\vee,\alpha_3^\vee, \alpha_4^\vee\}$ where now
$\alpha_1^\vee, \alpha_2^\vee$ are short, while $\alpha_3^\vee,
\alpha_4^\vee$ are long. Hence, Lusztig chooses (L) in \cite[4.10]{LuB}, 
while Alvis \cite[p.~6/22]{alvis}, Carter \cite[p.~414]{Ca2} and
{\sf CHEVIE} \cite{chevie}, \cite{jmich} choose~(AC).

Working with (AC) or (L) leads to different labels precisely for those 
irreducible characters that are permuted by the automorphism $\iota\colon W\rightarrow W$ in Remark~\ref{rem12}. Concretely, this means for example 
that the character $2_3$ has value $2$ on $s_1,s_2$ with respect to (L),
and value $0$ on $s_1,s_2$ with respect to (AC).
\end{rem}

%An explicit description of $\Irr(W)$ is given by Lusztig
%\cite[4.10]{LuB}, as follows.
%
%The trivial character is $1_1$, the sign character is $1_4$. By $1_2$ we
%denote the one-dimensional character which is $-1$ on $\tau,\tau\sigma$,
%and $+1$ on $d,a$; then $1_3=1_2\otimes \mbox{sgn}$. (Thus, we have
%$\mbox{sgn}=1_4=1_2\otimes 1_3$.) The character of the reflection
%representation is $4_2$; furthermore, $4_5=4_2\otimes \mbox{sgn}$,
%$4_3=4_2\otimes 1_2$ and $4_4=4_2\otimes 1_3$. The character $2_1$ is
%that of the characters of degree~$2$ which has value $0$ on $\tau,
%\tau\sigma$ and value~$2$ on $d,a$; moreover, $2_2=2_1\otimes \mbox{sgn}$.
%The character $2_3$ is that of the characters of degree~$2$ which has
%value $2$ on $\tau,\tau\sigma$ and value $0$ on $d,a$. Then $2_4=2_3\otimes 
%\mbox{sgn}$ and $4_1=2_1\otimes 2_3$. (Explicit constructions of $2_1$
%and $2_3$ will be given below.) The character $9_1$ can be described as
%the second symmetric power of $4_2$, from which we subtract the trivial
%character $1_1$; we have $9_4=9_1\otimes \mbox{sgn}$, $9_2=9_1\otimes 
%1_2$ and $9_3=9_1\otimes 1_3$. Furthermore, $8_1=4_2 \otimes 2_1$, $8_2=
%8_1\otimes \mbox{sgn}$, $8_3=4_2\otimes 2_3$ and $8_4=8_3 \otimes 
%\mbox{sgn}$. The character $6_2$ is the second exterior power of $4_2$
%and $6_1=6_2 \otimes 1_2=6_2\otimes 1_3$. Finally, $16_1=4_1 \otimes 4_2$
%and $12_1=6_1\otimes 2_1$.

%%%%%%%%%%%%%%%%%%%%%%%%%%%%%%%%%%%%%%%%%%%%%%%%%%%%%%%%%%%%%%%%%%%%%%%%%%%
\section{The Iwahori--Hecke algebra of type $F_4$} 
\label{sec1.5}

We shall also need the generic Iwahori--Hecke algebra $\cH$ associated 
with the Coxeter group $W$ of type $F_4$ as defined in Section~\ref{sec1}.
Let $K=\QQ(u,v)$ where $u,v$ are commuting indeterminates. 
Then $\cH$ is an associative $K$-algebra with basis $\{T_w\mid w\in W\}$. 
As an algebra, $\cH$ is generated by the basis elements $T_d,T_a,T_\tau, 
T_{\tau\sigma}$, subject to the above braid relations and the following
quadratic relations:
\begin{align*} 
T_d^2&=uT_1+(u-1)T_d,\qquad T_a^2=uT_1+(u-1)T_a,\\ T_\tau^2&= 
vT_1+(v-1)T_\tau,\qquad T_{\tau\sigma}^2=vT_1+(v-1)T_{\tau\sigma}.
\end{align*}
It is known that $\cH$ is abstractly isomorphic to the group algebra of 
$W$ over $K$ (an explicit isomorphism is provided by \cite{p115}); 
furthermore, the specialisation $(u,v)\mapsto (1,1)$ induces a 
bijection between the irreducible representations of $\cH$ and those of 
$W$. To give an example (which will also serve as a useful test case 
later on), we construct representations of $\cH$ which specialise to 
the representations of $W$ with characters $2_1$ and $2_3$. 

\begin{exmp} \label{ex12} By checking that the defininig relations for 
the above generators of $\cH$ hold, one obtains two-dimensional matrix 
representations $\sigma,\sigma'\colon \cH\rightarrow M_2(K)$ via the 
following assignments:
\begin{align*}
\sigma :& \; T_d\mapsto \left(\begin{array}{rr} u &  0\\ 0 & u
\end{array}\right), \;T_a\mapsto \left(\begin{array}{rr} u & 0 \\ 0 
& u\end{array}\right),\; T_{\tau}\mapsto \left(\begin{array}{rr}  
-1 & 1\\ 0 & v\end{array}\right), \;T_{\tau\sigma}\mapsto 
\left(\begin{array}{rr} v & 0 \\ v & -1\end{array}\right);\\
\sigma':& \; T_d\mapsto \left(\begin{array}{rr} u &  0\\ u & -1
\end{array}\right),\;T_a\mapsto \left(\begin{array}{rr} -1 & 1 \\ 
0 & u\end{array}\right),\; T_\tau\mapsto \left(\begin{array}{rr} v 
& 0 \\ 0 & v\end{array}\right),\; T_{\tau\sigma}\mapsto 
\left(\begin{array}{rr} v & 0 \\0 & v\end{array}\right).
\end{align*}
One easily sees that these are irreducible. We have 
\[\mbox{trace}(\sigma(T_d))=\mbox{trace}(\sigma(T_a))=2u\quad
\mbox{and}\quad \mbox{trace}(\sigma(T_\tau))=\mbox{trace}(\sigma
(T_{\tau\sigma}))=v-1.\]
Hence, if we specialise $(u,v)\mapsto (1,1)$, then the character of $\sigma$ 
becomes the character $2_1$ of $W$. Similarly, $\sigma'$ specialises to the 
character $2_3$ of $W$.

Now, we have the one-dimensional representation $\mbox{ind}\colon \cH
\rightarrow K$ such that $T_d,T_a\mapsto u$ and $T_\tau,T_{\tau\sigma}
\mapsto v$ (which specialises to the trivial character $1_1$ of $W$). 
Then the ``generic degree'' associated with an irreducible representation 
$\phi\colon \cH\rightarrow M_d(K)$ is defined by  
\[D_\phi:=\frac{d\sum_{w\in W} \mbox{ind}(T_w)}{\sum_{w\in W} \mbox{ind}
(T_w)^{-1}\mbox{trace}(\phi(T_w))\mbox{trace}(\phi(T_{w^{-1}}))} \in K
\qquad \mbox{(see \cite[p.~361]{Ca2})}.\]
All of these have been computed previously (see \cite[\S 13.5]{Ca2} and 
the references there) but, in the above two cases, we can just perform 
an explicit computation (using a computer) and obtain:
\begin{align*}
D_\sigma=D_{2_1} &=\frac{v^3(v+1)(uv^2+1)(u^2v^2+1)(u^3v^3+1)}{(u^3+1)
(u+v)(u^2+v)}, \\ D_{\sigma'}=D_{2_3} &= \mbox{ same formula but with the 
roles of $u,v$ exchanged}.
\end{align*}
(See also \cite[p.~450]{Ca2}.) Specialising $(u,v)$ to powers of $q$, 
where $q$ is a power of a prime, we obtain actual degrees of characters 
of certain Chevalley groups over $\FF_q$; see \cite[Theorem~10.11.5]{Ca2} and 
\cite[Corollary~8.7]{LuB}. 
\end{exmp}

%%%%%%%%%%%%%%%%%%%%%%%%%%%%%%%%%%%%%%%%%%%%%%%%%%%%%%%%%%%%%%%%%%%%%%%%%%%
\section{Harish-Chandra series of unipotent characters} \label{sec3}

Let $p$ be a prime and $G$ be a simple algebraic group over an algebraic 
closure of $\FF_p$. Let $q$ be a power of $p$ and $F\colon G \rightarrow 
G$ be a Frobenius map with respect to an $\FF_q$-rational structure. Let
$B\subseteq G$ be an $F$-stable Borel subgroup and $T\subseteq G$ be an
$F$-stable maximal torus contained in $B$. Let $\Phi$ be the root system
of $G$ with respect to $T$ and $\Pi\subseteq \Phi$ be the set of simple
roots determined by $B$. Then $F$ induces an action on $\Phi$ which 
preserves~$\Pi$. Let $G^F=G(\FF_q)$ and $\Uch(G^F)$ be the set of 
unipotent characters of $G^F$. We have a partition of $\Uch(G^F)$ into 
Harish-Chandra series, and there is a collection of bijections
\[ \{\cS\leftrightarrow \Irr(\cW_{\!\cS})\mid \cS \mbox{ Harish-Chandra
series of } \Uch(G^F)\}\]
where $\cW_{\!\cS}$ is a certain finite Coxeter group associated to each 
$\cS$. Following \cite[3.25]{cbms}, this is further specified as follows.
\begin{itemize}
\item There is a well-defined $F$-invariant subset $\Pi'\subseteq \Pi$ 
such that the simple reflections of $\cW_{\!\cS}$ are indexed by the 
set $\overline{\Pi}$ of $F$-orbits on $\Pi\setminus \Pi'$; furthermore, 
there is a certain ``parameter function'' $\lambda \colon \overline{\Pi} 
\rightarrow \{ q^n \mid n\in \ZZ_{\geq 1}\}$. 
\item The pair $(\cW_{\!\cS},\lambda)$ determines an Iwahori--Hecke algebra 
with parameters given by $\{\lambda(\bar{\alpha})\mid \bar{\alpha} \in 
\overline{\Pi}\}$ (see \cite[3.6]{cbms}); the degrees of the characters in 
$\cS$ are obtained from the generic degrees of that algebra using the 
formula in \cite[(3.26.1)]{cbms}.
\end{itemize}
The possibilities for $\Pi'\subseteq \Pi,\overline{\Pi},\lambda$ are 
explicitly listed in Table~II (p.~35) of \cite{cbms}. Let now $\cS\subseteq 
\Uch(G^F)$ be a Harish-Chandra series where the associated Coxeter group 
$\cW_{\!\cS}$ is of type $F_4$; this situation occurs for $G^F \in \{F_4(q), 
{^2\!E}_6(q),E_8(q)\}$.

\begin{exmp} \label{hc1} Let $G^F=F_4(q)$, where we fix the labelling of the
simple roots $\Pi\subseteq \Phi$ as in Section~\ref{sec1}:

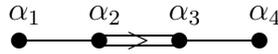
\vspace*{-4mm}
\begin{center}
\begin{picture}(220,20)
\put( 61,12){$\alpha_1$}
\put( 91,12){$\alpha_2$}
\put(121,12){$\alpha_3$}
\put(151,12){$\alpha_4$}
\put( 65,04){\circle*{6}}
\put( 95,04){\circle*{6}}
\put(125,04){\circle*{6}}
\put(155,04){\circle*{6}}
\put( 65,04){\line(1,0){30}}
\put( 95,02){\line(1,0){30}}
\put( 95,06){\line(1,0){30}}
\put(125,04){\line(1,0){30}}
\put(105,01){$>$}
\end{picture}
\end{center}
Let $\cS\subseteq \Uch(G^F)$ be the Harish-Chandra series containing the
trivial character of $G^F$. Then $\Pi'=\varnothing$, $\overline{\Pi}=\Pi$
and $\lambda(\alpha_i)=q$ for $i=1,2,3,4$. Thus, $\cW_{\!\cS}=W=\langle
s_1,s_2,s_3,s_4\rangle$ is the Weyl group of $G$ (where $s_i$ denotes 
the reflection with root $\alpha_i$); the characters in $\cS$ are in 
bijection with $\Irr(W)$.
As discussed in Section~\ref{sec1}, Lusztig chooses (L) in \cite[4.10]{LuB}, 
while Carter chooses (AC) in \cite[p.~414]{Ca2}, so one has to apply the 
automorphism in Remark~\ref{rem12} in order to pass from the notation 
concerning $\Uch(G^F)$ in Lusztig's book \cite[p.~371]{LuB} to Carter's
notation in Table~\ref{tabf4}.
\end{exmp}

\begin{rem} \label{hc1a} In the setting of Example~\ref{hc1}, let us
write $[\phi]\in \cS$ for the unipotent character corresponding to 
$\phi\in \Irr(W)$. If $(\phi',\phi'')$ is one of the pairs of characters 
that are permuted as in Remark~\ref{rem12}, then $[\phi'], [\phi'']$ 
have the same degree. Thus, as far as unipotent character degrees are 
concerned, it is irrelevant whether convention (AC) or convention (L) is
chosen. However, a difference occurs when we consider the partition of
the unipotent characters into ``families'' and the parametrisation of
the characters inside the various families in terms of finite sets
$\cM(\cG)$ (of certain pairs $(x,\sigma)$), as defined in
\cite[4.14]{LuB} (see also \cite[\S 13.6]{Ca2}). There are $11$ families,
two of which contain $4$ characters, one of which contains $21$ 
characters, and all the remaining families contain just one character. 
The $4$-element families contain $3$ characters of the form $[\phi]$ for 
$\phi\in \Irr(W)$. For example, one of them is given as follows:
\begin{center}
$\qquad$ \begin{tabular}{cc} \hline
\multicolumn{2}{c}{\mbox{Lusztig \cite[p.~371]{LuB}}} \\ \hline
$\left[4_2\right]$ & $(1,1)$ \\ $\left[2_1\right]$ & $(g_2,1)$ \\
$\left[2_3\right]$ & $(1,\varepsilon)$ \\ \hline
\end{tabular} $\qquad\qquad$ 
\begin{tabular}{cc} \hline
\multicolumn{2}{c}{\mbox{Carter \cite[p.~479]{Ca2}}} \\ \hline
$\left[ \phi_{4,1}\right]$ & $(1,1)$ \\ $\left[ \phi_{2,4}''\right]$ & 
$(g_2,1)$ \\ $\left[ \phi_{2,4}' \right]$ & $(1,\varepsilon)$ \\ \hline
\end{tabular} $\qquad\qquad$ 
\end{center}
Thus, the two labels $\phi_{2,4}''$ and $\phi_{2,4}'$ in Carter's table
should be exchanged, since Lusztig uses convention (L) and Carter uses
convention (AC). The same applies to the other $4$-element family, and
also to the $21$-element family.
\end{rem}

%is a family $\fF\subseteq \Uch(G^F)$ containing $21$ characters. The
%labelling for the $11$ characters in $\fF$ of the form $[\phi]$ 
%($\phi\in \Irr(W)$) are as follows (see \cite[p.~371]{LuB}): 
%\[ \renewcommand{\arraycolsep}{3pt} \begin{array}{cccccccccccc} \hline
%& (1,1) & (g_2',1) & (1,\lambda^1) & (g_2',
%\varepsilon') & ( 1,\lambda^2) & (g_2',\varepsilon'') & (g_4,1) & (g_2,
%\varepsilon'') & (g_3,1) & (1,\sigma) & (g_2,1) \\  \hline
%\mbox{(L)}: & [12_1] & [9_2] & [9_3] & [1_2] & [1_3] & [4_1] & [4_3] & 
%[4_4] & [6_1] & [6_2] & [16_1] \\ 
%\end{array}\]

\begin{exmp} \label{hc2} Let $G^F={^2\!E}_6(q)$. Then $\Phi$ is of type 
$E_6$ and $F$ induces an automorphism of order $2$ on $\Phi$. Let again 
$\cS\subseteq \Uch(G^F)$ be the Harish-Chandra series containing the 
trivial character of $G^F$. Then $\Pi'=\varnothing$ and $\overline{\Pi},
\lambda$ are given as follows.
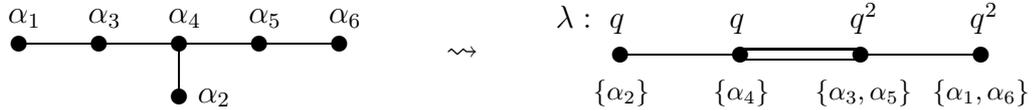
\begin{center}
\begin{picture}(395,40)
\put( 82, 0){$\alpha_2$}
\put( 11,30){$\alpha_1$}
\put( 41,30){$\alpha_3$}
\put( 71,30){$\alpha_4$}
\put(101,30){$\alpha_5$}
\put(131,30){$\alpha_6$}
\put(75, 2){\circle*{6}}
\put( 15,22){\circle*{6}}
\put( 45,22){\circle*{6}}
\put( 75,22){\circle*{6}}
\put(105,22){\circle*{6}}
\put(135,22){\circle*{6}}
\put(75,22){\line(0,-1){20}}
\put( 15,22){\line(1,0){30}}
\put( 45,22){\line(1,0){30}}
\put(75,22){\line(1,0){30}}
\put(105,22){\line(1,0){30}}
\put(175,16){$\leadsto$}

\put(216,28){$\lambda:$}
\put(236,28){$q$}
\put(281,28){$q$}
\put(326,28){$q^2$}
\put(371,28){$q^2$}
\put(230, 1){{\footnotesize $\{\alpha_2\}$}}
\put(275, 1){{\footnotesize $\{\alpha_4\}$}}
\put(313, 1){{\footnotesize $\{\alpha_3,\alpha_5\}$}}
\put(357, 1){{\footnotesize $\{\alpha_1,\alpha_6\}$}}
\put(240,18){\circle*{6}}
\put(285,18){\circle*{6}}
\put(330,18){\circle*{6}}
\put(375,18){\circle*{6}}
\put(240,18){\line(1,0){45}}
\put(285,16){\line(1,0){45}}
\put(285,20){\line(1,0){45}}
\put(330,18){\line(1,0){45}}
%\put(303,15){$>$}
\end{picture}
\end{center}
Again, $\cW_{\!\cS}$ is a Coxeter group of type $F_4$. When referring
to Kondo's article for the labelling of $\Irr(\cW_{\!\cS})$, one has 
to specify whether Kondo's generators $d,a$ correspond to $\{\alpha_2\},
\{\alpha_4\}$ or to $\{\alpha_3,\alpha_5\},\{\alpha_1,\alpha_6\}$. In 
other words, one has to specify whether $\lambda$ takes value $q$ on 
$d,a$ and value $q^2$ on $\tau,\tau\sigma$, or vice versa. This specification
can be reconstructed from the existing tables, as follows.

By \cite[Theorem~1.15]{Lu80}, there is a bijection between the unipotent
characters of the (untwisted) group $E_6(q)$ and those of ${^2\!E}_6(q)$;
if $\rho\in \Uch(E_6(q))$ corresponds to $\rho'\in \Uch({^2\!E}_6(q))$,
then the polynomial in $q$ which gives $\rho'(1)$ is obtained from the
polynomial which gives $\rho(1)$ by replacing $q$ by $-q$ (and adjusting
the sign). The bijection $\rho\leftrightarrow \rho'$ is defined by the tables
in \cite[1.10, 1.16]{Lu80}. In \cite{Lu80}, the characters of $\cS$ are
denoted $[\phi]$ for $\phi\in \Irr(\cW_{\!\cS})$. For example, using the 
information on character degrees in \cite[p.~363]{LuB}, we find that
\begin{align*}
\dim [2_1]&=\textstyle{\frac{1}{2}}q^3+ \mbox{higher powers of $q$},\\
\dim [2_3]&=q+ \mbox{higher powers of $q$}.
\end{align*}
On the other hand, as discussed above, the degrees of these two characters 
can also be obtained using the generic degrees of the Iwahori--Hecke 
algebra $\cH$ associated with $\cW_{\!\cS}$ and the parameter function 
$\lambda$.
Using the results of the computations at the end of 
Section~\ref{sec1.5}, and assuming $\lambda(d)=\lambda(a)=q$, 
$\lambda(\tau)=\lambda(\tau\sigma)=q^2$, we get the same result as above. 
(If we assume $\lambda(d)=\lambda(a)=q^2$, $\lambda(\tau)=\lambda(\tau\sigma)
=q$, then we get a different result.) Thus, Lusztig must have been
using the convention:

\begin{center}
\fbox{$\lambda(d)=\lambda(a)=q\quad\mbox{and}\quad \lambda(\tau)=
\lambda(\tau\sigma)=q^2$;}
\end{center}
and a comparison with the table in \cite[p.~481]{Ca2} shows that the same 
is also true for Carter.
Hence, both Lusztig \cite{LuB} and Carter \cite{Ca2} use the convention (AC) for the groups $G^F={^2\!E}_6(q)$.
\end{exmp}

\begin{exmp} \label{hc3} Let $G^F=E_8(q)$. There is a Harish-Chandra
series $\cS\subseteq \Uch(G^F)$ with corresponding $\Pi'\subseteq \Pi$, 
$\overline{\Pi},\lambda$ as follows (where $\Pi'$ is the subdiagram of
type $D_4$ indicated by open circles in the diagram for $E_8$):
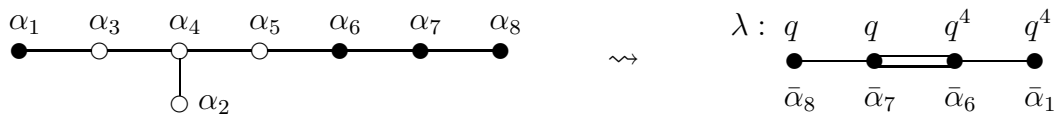
\begin{center}
\begin{picture}(405,40)
\put( 82, 0){$\alpha_2$}
\put( 11,30){$\alpha_1$}
\put( 41,30){$\alpha_3$}
\put( 71,30){$\alpha_4$}
\put(101,30){$\alpha_5$}
\put(131,30){$\alpha_6$}
\put(161,30){$\alpha_7$}
\put(191,30){$\alpha_8$}
\put(75, 2){\circle{6}}
\put( 15,22){\circle*{6}}
\put( 45,22){\circle{6}}
\put( 75,22){\circle{6}}
\put(105,22){\circle{6}}
\put(135,22){\circle*{6}}
\put(165,22){\circle*{6}}
\put(195,22){\circle*{6}}
\put(75,19){\line(0,-1){14}}
\put( 15,22){\line(1,0){27}}
\put( 48,22){\line(1,0){24}}
\put( 78,22){\line(1,0){24}}
\put(108,22){\line(1,0){89}}
\put(235,16){$\leadsto$}

\put(281,28){$\lambda:$}
\put(301,28){$q$}
\put(331,28){$q$}
\put(361,28){$q^4$}
\put(391,28){$q^4$}
\put(301, 1){$\bar{\alpha}_8$}
\put(331, 1){$\bar{\alpha}_7$}
\put(361, 1){$\bar{\alpha}_6$}
\put(391, 1){$\bar{\alpha}_1$}
\put(305,18){\circle*{6}}
\put(335,18){\circle*{6}}
\put(365,18){\circle*{6}}
\put(395,18){\circle*{6}}
\put(305,18){\line(1,0){30}}
\put(335,16){\line(1,0){30}}
\put(335,20){\line(1,0){30}}
\put(365,18){\line(1,0){30}}
%\put(345,15){$>$}
\end{picture}
\end{center}
Again, $\cW_{\!\cS}$ is a Coxeter group of type $F_4$. As in the previous
example, one has to specify whether $\lambda$ takes value $q$ on Kondo's
generators $d,a$ and value $q^4$ on $\tau,\tau\sigma$, or vice versa.
Actually, Lusztig specifies this in his book \cite[p.~361]{LuB}; it is the
same convention as in the previous example: 
\begin{center}
\fbox{$\lambda(d)=\lambda(a)=q\quad\mbox{and}\quad \lambda(\tau)=
\lambda(\tau\sigma)=q^4$;}
\end{center}
a comparison with the table in \cite[pp.~484--488]{Ca2} shows that the same 
is also true for Carter.
So both Lusztig \cite{LuB} and Carter \cite{Ca2} use the convention (AC) in the present situation.

One can also do a consistency check, as above. The characters in
$\cS$ are denoted by $D_4[\phi]$ for $\phi\in \Irr(\cW_{\!\cS})$. 
For example, according to the table in \cite[p.~366]{LuB}, the degrees 
of $D_4[2_1]$ and $D_4[2_3]$ should be
\begin{align*}
\dim D_4[2_1]&=\textstyle{\frac{1}{2}}q^{12}+ \mbox{higher powers of $q$},\\
\dim D_4[2_3]&=\textstyle{\frac{1}{2}}q^{4}+ \mbox{higher powers of $q$},
\end{align*}
On the other hand, we can also determine these degrees using the
generic degrees of the Iwahori--Hecke algebra associated with $\cW_{\!\cS},
\lambda$; we get the same result for the degrees precisely when we specify
the values of $\lambda$ as above. 
\end{exmp}

\begin{rem} \label{hc4} Let $\cW_{\!\cS}$ be the finite Coxeter group
associated with a Harish-Chandra series $\cS\subseteq \Uch(G^F)$. Then
\cite[Theorem~5.9]{Lu76} and \cite[(7.7)]{Lu76} do not only show how the 
type of $\cW_{\!\cS}$ and the parameter function $\lambda$ are determined 
from $\Pi'\subseteq \Pi$, but \cite[Theorem~5.9]{Lu76} also describes a 
natural root system $\overline{\Phi}$ for $\cW_{\!\cS}$. One can then check
the following:
\begin{itemize}
\item[(a)] In Example~\ref{hc2} (where $G^F={^2\!E}_6(q)$), the simple
roots in $\overline{\Phi}$ corresponding to the $F$-orbits $\{\alpha_2\}$,
$\{\alpha_4\}$ are long, and the other two are short.
\item[(b)] In Example~\ref{hc3} (where $G^F=E_8(q)$), the simple roots 
in $\overline{\Phi}$ corresponding to $\bar{\alpha}_8$, $\bar{\alpha}_7$
are long, and the other two are short.
\end{itemize}
Thus, as noted before, in both (a) and (b), the choice (AC) is made for
matching the simple reflections in $\cW_{\!\cS}$ (corresponding to the
roots in $\overline{\Pi}$) with Kondo's generators $d,a,\tau,\tau\sigma$. 
\end{rem}

%%%%%%%%%%%%%%%%%%%%%%%%%%%%%%%%%%%%%%%%%%%%%%%%%%%%%%%%%%%%%%%%%%%%%%%%%%%
\section{The (generalised) Springer correspondence} \label{genspr}

Let again $p$ be a prime and $G$ be a simple algebraic group over an 
algebraic closure of~$\FF_p$. Let $\cN_G$ be the set of all pairs $(C,\cE)$ 
where $C$ is a unipotent conjugacy class of~$G$ and $\cE$ is an irreducible 
$\overline{\QQ}_\ell$-local system on $C$, equivariant for the conjugation 
action of~$G$. (Here, $\ell$ is a prime $\neq p$.) Let $\cM_G$ be the
set of all triples $(L,C_0,\cE_0)$ (up to $G$-conjugacy) where $L$ is a
Levi subgroup of some parabolic subgroup of $G$ and $(C_0,\cE_0)\in\cN_L$
is ``cuspidal'' (see \cite[2.4, 6.2]{LuIC}). Using intersection cohomology
methods, Lusztig \cite[\S 6]{LuIC} has defined a natural surjective map 
$\cN_G\rightarrow \cM_G$. Following \cite{Lu20b}, the fibres of that map
will be called ``unipotent blocks''; they form a partition of~$\cN_G$. 
Furthermore, there is a collection of bijections
\[ \{\cI\leftrightarrow \Irr(\cW_{\!\cI})\mid \cI \mbox{ unipotent block of }
\cN_G\}\]
where $\cW_{\!\cI}$ is a certain finite Coxeter group associated to each 
$\cI$; in fact, by \cite[\S 9]{LuIC}, we have $\cW_{\!\cI}=N_G(L)/L$ where 
$(L,C_0,\cE_0)$ is the triple corresponding to $\cI$ under the map $\cN_G
\rightarrow\cM_G$. This is called the generalised Springer correspondence; 
note that all this does not require an $\FF_q$-rational structure on $G$. 
If $\cI_1$ is the unipotent block containing the pair $(\{1\},
\overline{\QQ}_\ell)$, then $\cW_{\!\cI_1}=W$ is the Weyl group of $G$ and 
the bijection $\cI_1 \leftrightarrow \Irr(W)$ is the correspondence defined 
earlier by Springer in the 1970s (up to tensoring by the sign character 
of $W$). For example, it is known that $\cI_1$ contains all pairs $(C,
\overline{\QQ}_\ell)$ where $C$ is a unipotent conjugacy class of~$G$
and $\overline{\QQ}_\ell$ stands for the trivial local system. We call
$\cI_1$ the ``principal unipotent block'' of $\cN_G$.

Now assume that $G$ is of exceptional type. The basic reference for the 
generalised Springer correspondence in this case is Spaltenstein's article
\cite{Spa1}, which completes and extends to small characteristics $p$ 
earlier work of Springer (for $G$ of type~$G_2$), Shoji ($F_4$) and 
Alvis--Lusztig ($E_6, E_7,E_8$); see the detailed references in \cite{Spa1}. 
According to \cite[\S 15]{LuIC}, there are three cases in which there are 
unipotent blocks $\cI\subseteq \cN_G$ such that~$\cW_{\!\cI}$ is of type 
$F_4$; these occur for $G$ of type 
\[ F_4 \;\mbox{ (any $p$)}, \qquad E_8 \;\mbox{ ($p=2$)}, \qquad  E_7
\;\mbox{ (simply connected and $p\neq 2$)}.\]
As in the previous section, there are certain issues concerning the
conventions used when referring to Kondo's labelling of $\Irr(\cW_{\!\cI})$.

\begin{exmp} \label{gspr1} Let $G$ be of type $F_4$, where we fix the 
labelling of the simple roots as in the previous section, that is, 
\begin{center}
\begin{picture}(220,20)
\put( 61,12){$\alpha_1$}
\put( 91,12){$\alpha_2$}
\put(121,12){$\alpha_3$}
\put(151,12){$\alpha_4$}
\put( 65,04){\circle*{6}}
\put( 95,04){\circle*{6}}
\put(125,04){\circle*{6}}
\put(155,04){\circle*{6}}
\put( 65,04){\line(1,0){30}}
\put( 95,02){\line(1,0){30}}
\put( 95,06){\line(1,0){30}}
\put(125,04){\line(1,0){30}}
\put(105,01){$>$}
\end{picture}
\end{center}
As above, let $\cI_1\subseteq \cN_G$ be the principal unipotent block.
Then $\cW_{\cI_1}=W=\langle s_1,s_2,s_3,s_4 \rangle$ is the Weyl group 
of $G$ (where $s_i$ is the reflection with root $\alpha_i$). The
correspondence $\cI_1\leftrightarrow \Irr(W)$ was determined by Shoji 
\cite{Sho}, with some conditions on the characteristic which were 
later removed by Spaltenstein \cite{Spa1}.  As before, one
has to specify which of Kondo's generators $d,a,\tau,\tau\sigma$ should
correspond to long roots, and which to short roots. Now Shoji explicitly 
describes the restrictions of the irreducible characters of~$W$ to a 
parabolic subgroup of type $C_3$. So one can deduce from \cite[4.3 and 
Table~4]{Sho} that Shoji uses the choice (AC) for 
the labelling of $\Irr(W)$. A comparison with the table in 
\cite[p.~330]{Spa1} shows that Spaltenstein uses the same choice. 
\end{exmp}

\begin{exmp} \label{gspr2} Let $G$ be of type $E_8$, where $p=2$. By 
\cite[15.3]{LuIC}, there is a unique unipotent block $\cI\subseteq \cN_G$ 
such that $\cW_{\!\cI}$ is a Coxeter group of type $F_4$. Under the
map $\cN_G\rightarrow \cM_G$, this block corresponds to a triple
$(L,C_0,\cE_0)$ where $L$ is of type $D_4$ (indicated by open 
circles in the diagram below). The generators of $\cW_{\!\cI}$ are 
indexed by the simple roots indicated by full circles in the diagram 
below (similarly to Example~\ref{hc3}):
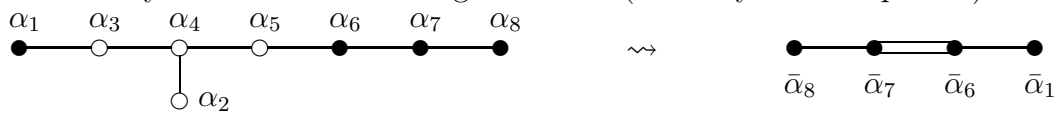
\begin{center}
\begin{picture}(405,35)
\put( 82, 0){$\alpha_2$}
\put( 11,30){$\alpha_1$}
\put( 41,30){$\alpha_3$}
\put( 71,30){$\alpha_4$}
\put(101,30){$\alpha_5$}
\put(131,30){$\alpha_6$}
\put(161,30){$\alpha_7$}
\put(191,30){$\alpha_8$}
\put(75, 2){\circle{6}}
\put( 15,22){\circle*{6}}
\put( 45,22){\circle{6}}
\put( 75,22){\circle{6}}
\put(105,22){\circle{6}}
\put(135,22){\circle*{6}}
\put(165,22){\circle*{6}}
\put(195,22){\circle*{6}}
\put(75,19){\line(0,-1){14}}
\put( 15,22){\line(1,0){27}}
\put( 48,22){\line(1,0){24}}
\put( 78,22){\line(1,0){24}}
\put(108,22){\line(1,0){89}}
\put(242,18){$\leadsto$}

\put(301, 5){$\bar{\alpha}_8$}
\put(331, 5){$\bar{\alpha}_7$}
\put(361, 5){$\bar{\alpha}_6$}
\put(391, 5){$\bar{\alpha}_1$}
\put(305,22){\circle*{6}}
\put(335,22){\circle*{6}}
\put(365,22){\circle*{6}}
\put(395,22){\circle*{6}}
\put(305,22){\line(1,0){30}}
\put(335,20){\line(1,0){30}}
\put(335,24){\line(1,0){30}}
\put(365,22){\line(1,0){30}}
%\put(345,15){$>$}
\end{picture}
\end{center}
Again, one has to match $\bar{\alpha}_8$, $\bar{\alpha}_7$, $\bar{\alpha}_6$,
$\bar{\alpha}_1$ to Kondo's generators $d,a,\tau,\tau\sigma$. Now 
Spaltenstein writes on \cite[p.~327]{Spa1} that it makes sense to assign
relative root lengths to $\bar{\alpha}_8$, $\bar{\alpha}_7$, $\bar{\alpha}_6$,
$\bar{\alpha}_1$. In this case, he declares $\bar{\alpha}_6$, $\bar{\alpha}_1$
to be long, and $\bar{\alpha}_8$, $\bar{\alpha}_7$ to be short. Since he
also refers to Alvis \cite{alvis} (with convention (AC)), we are led to the 
assumption that Spaltenstein uses the following matching with Kondo's 
generators:
\begin{center}
\fbox{$\quad \bar{\alpha}_1\leftrightarrow d, \qquad \bar{\alpha}_6
\leftrightarrow a, \qquad \bar{\alpha}_7\leftrightarrow \tau, \qquad 
\bar{\alpha}_8 \leftrightarrow \tau\sigma.\quad$}
\end{center}
Note that this is the matching opposite to that in Example~\ref{hc3}; 
note also that Spaltenstein's declaration of relative root lengths for
$\cW_{\!\cI}$ is opposite to that in Remark~\ref{hc4}(b). (These 
different declarations actually caused the ``delicate contradictions''
mentioned in the introduction, and prompted this note.)

That the above assumption is correct has been independently confirmed by 
the second author in \cite[\S 4.5]{Het3}, via a computation involving
characteristic functions of character sheaves and the intersection of
Bruhat cells with conjugacy classes. (If the opposite matching is used,
then those computations result in a contradiction; see 
\cite[Remark~4.5.33]{Het3}.) The fact that computations of that kind are
able to detect properties of the generalised Springer correspondence was 
subsequently used in order to resolve the last open question concerning 
the generalised Springer correspondence for $G$ of type $E_8$; see \cite{Het4}.
\end{exmp}

\begin{exmp} \label{gspr3} Let $G$ be simply connected of type $E_7$, where 
$p\neq 2$. By \cite[15.2]{LuIC}, there is a unique unipotent block $\cI
\subseteq \cN_G$ such that $\cW_{\!\cI}$ is a Coxeter group of type $F_4$. 
Under the map $\cN_G\rightarrow \cM_G$, this block corresponds to a triple
$(L,C_0,\cE_0)$ where $L$ is of type $A_1{+}A_1{+}A_1$ (indicated by open 
circles in the diagram below). The generators of $\cW_{\!\cI}$ are 
indexed by the simple roots indicated by full circles in the diagram below:
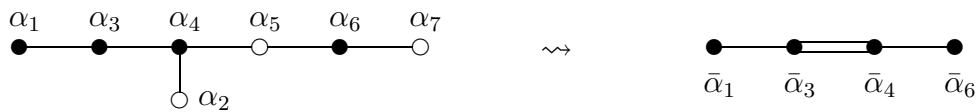
\begin{center}
\begin{picture}(405,40)
\put( 82, 0){$\alpha_2$}
\put( 11,30){$\alpha_1$}
\put( 41,30){$\alpha_3$}
\put( 71,30){$\alpha_4$}
\put(101,30){$\alpha_5$}
\put(131,30){$\alpha_6$}
\put(161,30){$\alpha_7$}
\put(75, 2){\circle{6}}
\put( 15,22){\circle*{6}}
\put( 45,22){\circle*{6}}
\put( 75,22){\circle*{6}}
\put(105,22){\circle{6}}
\put(135,22){\circle*{6}}
\put(165,22){\circle{6}}
\put(75,19){\line(0,-1){14}}
\put( 15,22){\line(1,0){27}}
\put( 48,22){\line(1,0){24}}
\put( 78,22){\line(1,0){24}}
\put(108,22){\line(1,0){54}}
\put(210,18){$\leadsto$}

%\put(251,28){$\lambda:$}
%\put(271,28){$1$}
%\put(301,28){$1$}
%\put(331,28){$4$}
%\put(361,28){$4$}
\put(271, 5){$\bar{\alpha}_1$}
\put(301, 5){$\bar{\alpha}_3$}
\put(331, 5){$\bar{\alpha}_4$}
\put(361, 5){$\bar{\alpha}_6$}
\put(275,22){\circle*{6}}
\put(305,22){\circle*{6}}
\put(335,22){\circle*{6}}
\put(365,22){\circle*{6}}
\put(275,22){\line(1,0){30}}
\put(305,20){\line(1,0){30}}
\put(305,24){\line(1,0){30}}
\put(335,22){\line(1,0){30}}
%\put(315,15){$>$}
\end{picture}
\end{center}
Similarly to the previous example, Spaltenstein \cite[p.~327]{Spa1} declares
$\bar{\alpha}_4$, $\bar{\alpha}_6$ to be long, and $\bar{\alpha}_1$, 
$\bar{\alpha}_3$ to be short. (Note again that this declaration of relative 
root lengths is opposite to that provided by \cite[Theorem~5.9]{Lu76}.) Again,
we are led to assume that Spaltenstein uses the following matching 
with Kondo's generators:
\begin{center}
\fbox{$\quad \bar{\alpha}_6\leftrightarrow d, \qquad \bar{\alpha}_4
\leftrightarrow a, \qquad \bar{\alpha}_3\leftrightarrow \tau, \qquad 
\bar{\alpha}_1 \leftrightarrow \tau\sigma.\quad$}
\end{center}
This is confirmed by the computations in the proof of Lemma~\ref{esp}
below, in which also one error in Spaltenstein's table will be corrected. 
We came across that error via the following remark.
\end{exmp}

%A consistency check is as follows. Consider the unipotent class $C$
%denoted by $A_5''$ in Spaltenstein's table \cite[p.~331/332]{Spa1}. 
%For $u\in C$, we have $A_G(u):=C_G(u)/C_G^\circ(u)\cong \ZZ/2\ZZ$.
%Furthermore, the pair $(C,\cE)\in \cN_G$, where $\cE$ corresponds to 
%the non-trivial character of $A_G(u)$, belongs to $\cI$ and we have 
%$(C,\cE)\leftrightarrow 8_4$. By the definition of the notation,
%we have $C_1\subseteq C$ where $C_1$ is the class of regular unipotent 
%elements in the Levi subgroup~$M\subseteq G$ of type $A_5$ corresponding 
%to the simple roots $\alpha_2, \alpha_4,\alpha_5, \alpha_6,\alpha_7$; 
%furthermore, $A_M(u_1):=C_M(u_1)/C_M^\circ(u_1)\cong \ZZ/2\ZZ$ for 
%$u_1\in C_1$. Let $\cI_1$ be the unipotent block of $\cN_M$ containing 
%the pair $(C_1,\cE_1)$ where $\cE_1$ corresponds to the non-trivial 
%character of $A_M(u_1)$; here, $\cW_{\!\cI_1}$ is the parabolic subgroup 
%of type $A_2$ corresponding to the roots $\bar{\alpha}_4$, 
%$\bar{\alpha}_6$ of $\cW_{\!\cI}$. Under the bijection $\cI_1
%\leftrightarrow \Irr(\cW_{\!\cI_1})$, the corresponding character of 
%$\cW_{\!\cI_1}$ is the trivial character. By the multiplicity formula 
%\cite[1.2(II)]{Spa1} and the discussion in \cite[1.4]{Spa1}, the trivial 
%character of $\cW_{\!\cI_1}$ occurs with multiplicity~$1$ in the 
%restriction of the character $8_4$ of $\cW_{\!\cI}$. Using {\sf CHEVIE},
%we find that this is the case only when $d,a$ correspond to long roots. 

\begin{rem} \label{lurem} Let $\cI\subseteq \cN_G$ be a unipotent block. 
Using the bijection $\cI\leftrightarrow \Irr(\cW_{\!\cI})$, we introduce 
an equivalence relation $\sim$ on $\Irr (\cW_{\!\cI})$ as follows. Let 
$\chi,\chi'\in \Irr(\cW_{\!\cI})$ and assume that $(C,\cE)\leftrightarrow 
\chi$ and $(C',\cE')\leftrightarrow \chi'$ under the above bijection; then 
we write $\chi\sim \chi'$ if $C=C'$. In \cite{Lu20a}, Lusztig shows that the 
equivalence classes under $\sim$ for the principal unipotent block $\cI_1
\subseteq \cN_G$ can be recovered in a purely algebraic way, using operations 
with characters of Weyl groups. In \cite[Conj.~6.11]{Lu20b}, there is a 
conjecture in a similar spirit for arbitrary unipotent blocks of $\cN_G$. 
Using {\sf CHEVIE} \cite{jmich}, we computed the algebraic version of 
$\sim$ for the unipotent block $\cI\subseteq \cN_G$ in Example~\ref{gspr3}, 
where $G$ is of type $E_7$; and we found an inconsistency with Spaltenstein's
table \cite{Spa1}.
\end{rem}

%%%%%%%%%%%%%%%%%%%%%%%%%%%%%%%%%%%%%%%%%%%%%%%%%%%%%%%%%%%%%%%%%%%%%%%%%%%
\section{Computing unipotent blocks}

In this section we briefly discuss the explicit computation of unipotent 
blocks and the generalised Springer correspondence. In particular, this
will allow us to resolve the inconsistency in Spaltenstein's table 
\cite{Spa1} mentioned in Remark~\ref{lurem}.

The main tool is the ``multiplicity formula'' in \cite[1.2(II)]{Spa1} 
(which is a reformulation of a formula that originally appeared in 
\cite[\S 8]{LuIC}). It relates the collection of bijections
\[ \{\cI\leftrightarrow \Irr(\cW_{\!\cI})\mid \cI \mbox{ unipotent block of }
\cN_G\}\]
to the analogous collection of bijections
\[ \{\cI'\leftrightarrow \Irr(\cW_{\!\cI'})\mid \cI' \mbox{ unipotent
block of } \cN_{M}\}\]
where $M$ is a Levi subgroup of some parabolic subgroup of $G$. Recall 
that we have surjective maps $\cN_G\rightarrow \cM_G$ and $\cN_M 
\rightarrow \cM_M$. Let us fix unipotent blocks $\cI\subseteq \cN_G$ 
and $\cI'\subseteq \cN_M$ which correspond to the same triple in $\cM_M$.
(Note that $\cM_M\subseteq \cM_G$.) The latter condition implies that 
$\cW_{\!\cI'}$ can be naturally regarded as a parabolic subgroup 
of~$\cW_{\!\cI}$. To state the multiplicity formula, we need to 
introduce some further notation. 

Let $(C,\cE) \in \cI$ and $\rho\in \Irr(\cW_{\!\cI})$ correspond to
$(C,\cE)$ under the above bijection. Furthermore, we set $A_G(u):=C_G(u)/
C_G^\circ(u)$ where $u\in C$; then $\cE$ corresponds to an irreducible 
character $\phi\in \Irr(A_G(u))$ (see \cite[\S 0]{LuIC}). Let also 
$(C',\cE')\in \cI'$ and $\rho'\in \Irr(\cW_{\!\cI'})$ correspond to 
$(C',\cE')$ under the above bijection. Furthermore, we set $A_M(u'):=
C_M(u')/C_M^\circ(u')$ where $u'\in C'$; then $\cE'$ corresponds to an 
irreducible character $\phi'\in \Irr(A_M(u'))$. In this setting, the 
multiplicity formula states:
\begin{align*}
& \mbox{multiplicity of $\rho'\in \Irr(\cW_{\!\cI'})$ in the restriction 
of $\rho\in \Irr(\cW_{\!\cI})$ to $\cW_{\!\cI'}$}\\
=\; & \mbox{multiplicity of $\phi\otimes \overline{\phi}'\in 
\Irr(A_G(u)\times A_M(u'))$ in $\varepsilon_{u,u'}$},
\end{align*}
where $\varepsilon_{u,u'}$ is the character of a permutation representation 
of the direct product $A_G(u)\times A_M(u')$ on a certain finite set 
$X_{u,u'}$ (defined in \cite[1.2]{Spa1}). 

\begin{exmp} \label{xuu1a} In the above setting, let $\cI_1$ be the principal
unipotent block of $\cN_G$ and $\cI_1'$ be the principal unipotent block of 
$\cN_M$. Assume that the bijections $\cI_1\leftrightarrow\Irr(\cW_{\!\cI_1})$ 
and $\cI_1'\leftrightarrow \Irr(\cW_{\!\cI_1'})$ (that is, the ordinary 
Springer correspondences) are already known. Now, as remarked earlier, we 
have $(C,\overline{\QQ}_\ell) \in \cI_1$ and $(C',\overline{\QQ}_\ell)\in 
\cI_1'$. The corresponding characters of $A_G(u)$ and $A_M(u')$ are the 
trivial characters. So, since $\varepsilon_{u,u'}$ is the character of a 
permutation representation, the right hand side of the multiplicity formula 
will be a positive integer if $X_{u,u'}\neq \varnothing$, and $0$ otherwise.
Via the left hand side of the multiplicity formula, this condition can be 
checked using a computation with the characters of $\cW_{\!\cI_1}$ and of
$\cW_{\!\cI_1'}$. Assume now that we have a situation where $X_{u,u'}=
\varnothing$ and $\cI\neq \cI_1$, $\cI'\neq \cI_1'$. Then the right 
hand side of the multiplicity formula for the pairs $(C,\cE)\in \cI$ and
$(C',\cE')\in \cI'$ will still be $0$ (because $X_{u,u'}=\varnothing$) and 
so the multiplicity of $\rho'\in \Irr(\cW_{\!\cI'})$ in the restriction 
of $\rho\in \Irr(\cW_{\!\cI})$ to $\cW_{\!\cI'}$ must be $0$. If the
bijection $\cI'\leftrightarrow \Irr(\cW_{\!\cI'})$ is already known, then
this rules out many cases on the level of the bijection $\cI\leftrightarrow 
\Irr(\cW_{\!\cI})$~---~as pointed out in \cite[4.3]{Spa1}.
\end{exmp}

\begin{lem} \label{esp} Let $G$ and $\cI\subseteq \cN_G$ be as in 
Example~\ref{gspr3}, where $G$ is simply connected of type $E_7$ (with 
$p\neq 2$) and $\cW_{\!\cI}$ is of type $F_4$. 
\begin{itemize}
\item[(a)] In the table \cite[p.~331/332]{Spa1}, the characters $\chi_{2,3}$ 
and $\chi_{8,3}$ of $\cW_{\!\cI}$ should be interchanged in the column 
labelled $3A_1'' \;(p\neq 2)$.
\item[(b)] With the adjustment in (a), the statements in
\cite[Conj.~6.11]{Lu20b} are correct for~$\cI$.
\end{itemize}
\end{lem}

(Here, we use Spaltenstein's notation for the irreducible
characters of $\cW_{\!\cI}$.)

\begin{proof} (a) We only sketch this. First of all, by \cite[15.2]{LuIC},
$\cN_G$ is the union of $\cI_1, \cI$ and one block consisting of a
cuspidal pair when $p\neq 2,3$ (where, as usual, $\cI_1$ denotes the 
principal unipotent block); if $p=3$, then $\cN_G$ is the union of
$\cI_1,\cI,\cI_2$ and three blocks consisting of cuspidal pairs, 
where $\cI_2$ is a block with $\cW_{\!\cI_2}$ of type $A_1$. We shall 
assume that $\cI_1$ has been determined. By \cite[5.5]{Spa1}, the block 
$\cI_2$ is also easily determined since $\cW_{\!\cI_2}$ is of type $A_1$;
the cuspidal pairs are listed in the last table of \cite[p.~337]{Spa1}. 
Thus, we can assume that the subset $\cI\subseteq \cN_G$ is known. 

Now let $C$ be the unipotent class of $G$ denoted by $D_5{+}A_1$, where 
$\dim \fB_u=6$ and $A_G(u)\cong \ZZ/2\ZZ$ for $u\in C$. Consider the pair
$(C,\cE)\in \cN_G$, where $\cE$ corresponds to the non-trivial character 
of $A_G(u)$. We have $(C,\cE)\in \cI$. According to Spaltenstein 
\cite[p.~174]{Spa0}, $C$ is obtained by the process of induction 
\cite{LuSp} from the trivial class of a Levi of type $A_2{+}A_2$.
By the transitivity of induction, it will also be induced from a class 
$C'$ of a Levi $M$ of type $A_5$ containing $A_2{+}A_2$. We choose $M$ 
such that it contains the Levi~$L$ of type $A_1{+}A_1{+}A_1$ (indicated 
by open circles in Example~\ref{gspr3}). Using weighted Dynkin diagrams and 
\cite[Prop.~1.9(b)]{LuSp}, we see that~$C'$ is the unique unipotent class 
of $M$ with $\dim \fB_{u'}=6$ und $A_M(u')\cong\ZZ/2\ZZ$ ($u'\in C'$). 

Consider the pair $(C',\cE')\in \cN_M$ where $\cE'$ corresponds to the 
non-trivial character of $A_M(u')$. The following information on $(C',
\cE')$ is obtained via the {\sf CHEVIE} function {\tt UnipotentClasses}. 
The pair $(C',\cE')$ belongs to the unipotent block $\cI'\subseteq \cN_{M}$ 
which corresponds to the same triple in $\cM_M\subseteq \cM_G$ as $\cI$.
Furthermore, the parabolic subgroup $\cW_{\!\cI'}\subseteq \cW_{\!\cI}$ 
is of type $A_2$, with simple roots corresponding to $\bar{\alpha}_4, 
\bar{\alpha}_6$ (long roots, as in Example~\ref{gspr3}). Under the 
bijection $\cI' \leftrightarrow \Irr(\cW_{\!\cI'})$, the corresponding
character is the sign character of $\cW_{\!\cI'}$. 

Since $C$ is obtained from $C'$ by induction, we are in the set-up
of \cite[1.3]{Spa1}. This shows that there are subgroups $H\subseteq N 
\subseteq A_G(u)$ such that $H$ is normal in $N$ with $N/H\cong A_M(u')$;
furthermore, the sets $X_{u,u'}$ and $A_G(u)/H$ are isomorphic as
sets with $(A_G(u)\times A_M(u'))$-actions. Since $A_G(u)$ and $A_M(u')$ 
are both isomorphic to $\ZZ/2\ZZ$, we must have $H=\{1\}$. Thus,
$X_{u,u'}$ and $A_G(u)$ are isomorphic as sets with $(A_G(u) \times 
A_M(u'))$-actions. Consequently, we find that the right hand side of the 
multiplicity formula evaluates to~$1$. Working out the left hand side of 
that formula (using the {\sf CHEVIE} function {\tt InductionTable}), we 
see that only the following characters of $\cW_{\!\cI}$ can correspond to
$C$: $\chi_{1,3}$, $\chi_{1,4}$, $\chi_{6,1}$, $\chi_{6,2}$, $\chi_{8,3}$, 
$\chi_{8,4}$. Using also induction of suitable classes from Levi subgroups 
$M$ of type $A_5{+}A_1$ and $D_6$, one rules out further characters from 
the above list until the only remaining possibility is $\chi_{8,3}$.
Thus, under the bijection $\cI\leftrightarrow \Irr(\cW_{\!\cI})$, the
pair $(C,\cE)$ must correspond to the character $\chi_{8,3}$.

For the remaining characters of $\cW_{\!\cI}$, one proceeds as follows.
First one considers all unipotent classes $C$ such that $|A_G(u)|=2$ for
$u\in C$. Using similar arguments as above, plus the highly efficient 
method in Example~\ref{xuu1a}, one checks that the entries for these classes 
in Spaltenstein's table are correct. It then remains 
to consider the classes $C$ denoted 
\[ D_6{+}A_1, \quad D_6(a_1){+}A_1,  \quad D_6(a_2){+}A_1, \quad 
D_4(a_1){+}A_1,\]
where only the following characters of $\cW_{\!\cI}$ can correspond to any
of these classes: 
\[\chi_{1,2},\quad \chi_{2,3}, \quad \chi_{2,4}, \quad \chi_{4,2},\quad 
\chi_{6,2}, \quad \chi_{8,1},\quad \chi_{9,4},\quad \chi_{12}.\]
We now use the method in Example~\ref{xuu1a}, with respect to a Levi $M
\subseteq G$ of type $D_6$, in order to rule out a number of possibilities. 
Firstly, we find that only $\chi_{1,2}, \chi_{8,1}$ from the above list 
can correspond to $D_6{+}A_1$ (in accordance with Spaltenstein's 
table); next, one finds that only $\chi_{1,2},\chi_{2,3},\chi_{4,2}$ can 
correspond to $D_6(a_1){+}A_1$. So the conclusion is that $\chi_{2,3},
\chi_{4,2}$ must correspond to $D_6(a_1){+}A_1$. Similarly, one finds 
that only $\chi_{2,4},\chi_{9,4}$ can correspond to $D_4(a_1){+}A_1$ and
that only $\chi_{2,4},\chi_{6,2},\chi_{12}$ can correspond to $D_6(a_2)
{+}A_1$. Hence, again, the conclusion is that $\chi_{6,2},\chi_{12}$ must 
correspond to $D_6(a_2){+}A_1$. Thus, all entries in Spaltenstein's table 
for $E_7$ are found to be correct, with the exception that $\chi_{2,3}$ 
and $\chi_{8,3}$ need to be exchanged. 

%For $D_6(a_1){+}A_1$ method in \cite[4.3]{Spa1}:
%Use $M$ of type $D_6$ rel. group of type $C_3$. $u' in M$ type $4431$
%mit Springer $2.22$. Skalar product with $189_5$ is 0, so  set $X_{u,u'}$ 
%is empty. But scalar product of $8_3$ in rel group non-zero, contradiction.
%
%(Or with regular unip class of $M$, rules out $8_1$.)

(b) This simply follows by inspection of the results of
the computation.
\end{proof}

\medskip
\noindent {\bf Acknowledgements}. We thank George Lusztig and Toshiaki
Shoji for useful comments. We also thank Gunter Malle for comments on an 
earlier version.  This article is a contribution to Project-ID 
286237555~--~TRR 195~--~by the Deutsche Forschungsgemeinschaft (DFG, 
German Research Foundation).

%%%%%%%%%%%%%%%%%%%%%%%%%%%%%%%%%%%%%%%%%%%%%%%%%%%%%%%%%%%%%%%%%%%%%%%%%%%

\end{document}